\documentclass[a5paper,12pt]{article}

\usepackage[a5paper,includehead,margin=2mm,headheight=13.6pt]{geometry}

\makeatletter
\@ifpackageloaded{amssymb}{}{\usepackage[bitstream-charter]{mathdesign}}
\makeatother

\usepackage{amsmath}

\newcommand{\secformat}[1]{\centering{\normalfont\normalsize{#1}}}

\usepackage{fancyhdr}

 {
	\fancyhf{}
	\fancyhead[L]{\small\thepage}
	\fancyhead[C]{\small\textsc{Frobenius}: A generalization of \textsc{Sylow}'s theorem.}
	\fancyhead[R]{\small[\firstmark--\botmark]}
}
\fancypagestyle{plain}{%
	\fancyhf{}
	\fancyhead[L]{\small\thepage}
	\fancyhead[R]{\small[\firstmark--\botmark]}
}

\usepackage{indentfirst}

\let\fr\mathfrak
\newcommand{\CA}{\fr{A}}
\newcommand{\CB}{\fr{B}}

\newcommand{\CD}{\fr{D}}
\newcommand{\CG}{\fr{G}}
\newcommand{\CH}{\fr{H}}
\newcommand{\CK}{\fr{K}}
\newcommand{\CP}{\fr{P}}
\newcommand{\CQ}{\fr{Q}}
\newcommand{\CR}{\fr{R}}

\usepackage{color}
\definecolor{darkblue}{rgb}{0, 0, 0.5}

\usepackage[%
		pdftitle={A generalization of Sylow's theorem (G. Frobenius)},
		pdfauthor={G. Frobenius},
		bookmarksopen,
		colorlinks, linkcolor=darkblue, urlcolor=darkblue, citecolor=darkblue,
		pagebackref
]{hyperref}

\newcounter{origpagecounter}{}
\newcommand{\origpagefirst}[1]{\smash{\setcounter{origpagecounter}{#1}\mark{\arabic{origpagecounter}}}}
\newcommand{\origpagebreak}{\mark{\arabic{origpagecounter}}\addtocounter{origpagecounter}{1}\mark{\arabic{origpagecounter}}}

\begin{document}

\pagenumbering{alph}
\thispagestyle{empty}

\begin{center}
	\begin{minipage}{0.7\textwidth}
		A translation of
		\begin{quote}
			``Verallgemeinerung des Sylow'schen Satzes''
			by
			F.~G.~Frobenius
			\\
			Sitzungsberichte der
			K\"onigl.~Preu\ss.~Akad.~der Wissenschaften
			zu Berlin,
			1895 (II),
			981--993.
		\end{quote}
	\end{minipage}
\end{center}

\clearpage

\pagenumbering{arabic}
\thispagestyle{plain}
\setcounter{page}{1}
\origpagefirst{981}

\begin{center}
	{\Large
		{A generalization of \textsc{Sylow}'s theorem.}
		\\
	}
	{\large
		{By \textsc{G.~Frobenius}.}
		\\
		---------
	}
\end{center}



\noindent
Every finite group whose order
is divisible by the prime $p$
contains elements of order $p$.
(\textsc{Cauchy},
\emph{M\'emoire sur les arrangements que l'on peut former avec des lettres donn\'ees}.
Exercises d'analyse et de physique Math\'ematique,
Vol.~III,
\S.~XII,
p.~250.)
Their number is,
as I will show here,
always a number of the form
$(p - 1) (n p + 1)$.
From that theorem,
\textsc{Sylow} deduced the more general one,
that a group
whose order $h$ is divisible by $p^\kappa$,
must contain subgroups of order $p^\kappa$.
(\emph{Th\'eor\`emes sur les groupes de substitutions},
Math.~Ann.,~Vol.~V.)
I gave a simple proof thereof
in my work
\emph{Neuer Beweis des \textsc{Sylow}'schen Satzes},
\textsc{Crelle}'s Journal,
Vol.~100.
The number of those subgroups must,
as I will show here,
always be $\equiv 1 \pmod{p}$.
If $p^\lambda$ is the highest power of $p$
contained in $h$,
then \textsc{Sylow} proved this theorem
only for the case that $\kappa = \lambda$.
Then any two groups of order $p^\lambda$
contained in $\CH$
are conjugate,
and their number $n p + 1$
is a divisor of $h$,
while for $\kappa < \lambda$
this does not hold in general.
I obtain the stated results
in a new way
from a theorem of group theory
that appears to be unnoticed thus far:

\emph{ %
In a group of order $h$,
the number of elements
whose order divides $g$
is divisible by
the greatest common divisor
of $g$ and $h$.
}

\subsubsection*{\secformat{\S.~1.}}

If $p$ is a prime number
then any group $\CP$ of order $p^\lambda$
has
a series of invariant subgroups (chief series)
$\CP_1, \CP_2, \ldots, \CP_{\lambda-1}$
of orders
$p, p^2, \ldots, p^{\lambda-1}$,
each of which is contained in the subsequent one.
\textsc{Sylow} (loc.~cit., p.~588) derives this result
from the theorem:

I.
\phantomsection\label{t:1-1}
\emph{ %
Every group of order $p^\lambda$ contains
an invariant element of order $p$.
}

An invariant element of a group $\CH$
is an element of $\CH$
that
is permutable with every element of $\CH$.
\origpagebreak
If $\CP$ contains
the invariant element $P$ of order $p$
then
the powers
of $P$ form
an invariant subgroup $\CP_1$ of $\CP$
whose order is $p$.
Likewise,
$\CP / \CP_1$ has
an invariant subgroup $\CP_2 / \CP_1$ of order $p$
hence
$\CP$ has an invariant subgroup $\CP_2$ of order $p^2$
which contains $\CP_1$, etc.
In my work
\emph{\"Uber die Congruenz nach einem aus zwei endlichen Gruppen gebildeten Doppelmodul},
\textsc{Crelle}'s Journal, Vol.~101
(\S.~3, IV),
I complemented
that theorem with the following remark:

II.
\phantomsection\label{t:1-2}
\emph{ %
Every group of order $p^{\lambda-1}$
contained in a group of order $p^\lambda$
is an invariant subgroup.
}

Other proofs for this
I developed
in my work
\emph{\"Uber endliche Gruppen},
Sit\-zungs\-be\-richte 1895
(\S.~2, III, IV, V; \S.~4, II).
This can be obtained from
\hyperref[t:1-1]{Theorem I}
in the following way:
Let
$\CH$ be a group of order $p^\lambda$,
$\CG$ a subgroup of order $p^{\lambda-1}$,
$P$ an invariant element of $\CH$
whose order is $p$,
and
$\CP$ the group of the powers of
$P$. 
If $\CG$ is divisible by $\CP$
then
$\CG / \CP$
is an invariant subgroup of $\CH / \CP$
because
on can assume
\hyperref[t:1-2]{Theorem II}
as proven
for groups whose order is smaller than $p^\lambda$.
Thus $\CG$ is an invariant subgroup of $\CH$.
If $\CG$ is not divisible by $\CP$
then $\CH = \CG \CP$,
meaning that
every element of $\CH$
can be brought into the form $H = G P$,
where $G$ is an element of $\CG$.
Now,
$G$ is permutable with $\CG$
and
$P$ even with every element of $\CG$.
Hence
also $H$ is permutable with $\CG$.

The theorem mentioned
at the onset
lends itself
to completion
in a different direction:

III.
\phantomsection\label{t:1-3}
\emph{ %
Every invariant subgroup of order $p$
of a group of order $p^\lambda$
consists
of powers of an invariant element.
}

Let $\CH$ be a group of order $p^\lambda$,
$\CP$ an invariant subgroup of order $p$.
If $Q$ is any element of $\CH$
and $q = p^\kappa$ is its order,
then
the powers of $Q$ form
a group $\CQ$ contained in $\CH$ of order $q$.
If $\CP$ is a divisor of $\CQ$
then every element $P$ of $\CP$
is a power of $Q$,
hence permutable with $Q$.
If $\CP$ is not a divisor of $\CQ$
then $\CP$ and $\CQ$
are relatively prime.
$\CP$ is permutable
with every element of $\CH$
and therefore
with every element of $\CQ$.
Thence $\CP \CQ$ is a group of order $p^{\kappa + 1}$
and
$\CP$ is an invariant subgroup of it.
But by
\hyperref[t:1-2]{Theorem II},
$\CQ$ is one also.
Therefore
$P$ and $Q$ are permutable
in view of the Theorem:

IV.
\phantomsection\label{t:1-4}
\emph{ %
If each of the relatively prime
groups $\CA$ and $\CB$
is permutable with every element of the other,
then
every element of $\CA$ is permutable
with
every element of $\CB$.
}

\origpagebreak

Indeed,
if $A$ is an element of $\CA$
and $B$ is an element of $\CB$,
then
the element
\[
	A ( B A^{-1} B^{-1} ) = (A B A^{-1}) B^{-1}
\]
is
contained in both $\CA$ and $\CB$,
and is therefore
the
principal element $E$.

I want to prove
\hyperref[t:1-3]{Theorem III}
also in a second way:
If $Q^{-1} P Q = P^a$
then
$Q^{-q} P Q^q = P^{a^q}$.
Hence if $Q^q = E$ then $a^q \equiv 1 \pmod{p}$.
Now $a^{p-1} \equiv 1 \pmod{p}$,
hence
as $q$ and $p-1$ are relatively prime,
also $a \equiv 1 \pmod{p}$
and therewith
$P Q = Q P$.

Thirdly and finally,
the Theorem follows from the more general Theorem:

V.
\phantomsection\label{t:1-5}
\emph{ %
Every invariant subgroup of a group $\CH$ of order $p^\lambda$
contains an invariant element of $\CH$
whose order is $p$.
}

Partition the elements of $\CH$
into classes of conjugate elements
(conjugate with respect to $\CH$).
If a class consists of a single element,
then it is an invariant one,
and conversely
every invariant element of $\CH$ forms a class by itself.
Let $\CG$ be an invariant subgroup of $\CH$
and $p^\kappa$ its order.
If the group $\CG$ contains an element of a class
then it contains all its elements.
Select an element $G_1, G_2, \ldots, G_n$
from each of the $n$ classes contained in $\CG$. 
If the elements of $\CH$ permutable with $G_\nu$
form a group of order $p^{\lambda_\nu}$,
then the number of elements
of $\CH$ conjugate to $G_\nu$,
i.e.~the number of elements
in the class represented by $G_\nu$, 
equals $p^{\lambda - \lambda_\nu}$
(\textsc{Crelle}'s Journal, Vol.~100, p.~181).
Thence
\[
	p^\kappa
	=
	p^{\lambda - \lambda_1}
	+
	p^{\lambda - \lambda_2}
	+
	\cdots
	+
	p^{\lambda - \lambda_n}
	.
\]

If $G_1$ is the principal element $E$
then $\lambda = \lambda_1$.
Therefore
not all the last $n-1$ terms
on the right hand side of this equation
can be divisible by $p$.
There must exist therefore another index $\nu > 1$
for which $\lambda_\nu = \lambda$ holds.
Then $G_\nu$ is an invariant element of $\CH$
whose order is $p^\mu > 1$,
and
the $p^{\mu-1}$-th power of $G_\nu$
is an invariant element of $\CH$ of order $p$
that is contained in $\CG$.

\subsubsection*{\secformat{\S.~2.}}

I.
\phantomsection\label{t:2-1}
\emph{ %
If $a$ and $b$ are relative primes,
then any element of order $a \, b$
can always, and in a unique way,
be written
as a product of two elements
whose orders are $a$ and $b$
and
which are permutable with each other.
}

If $A$ and $B$ are two permutable elements
whose orders $a$ and $b$ are relative primes,
then $A B = C$ has the order $a b$.
Conversely, let $C$ be any element of order $a b$.
Determining the integer numbers $x$ and $y$
such that $a x + b y = 1$
and
setting $a x = \beta$, $b y = \alpha$,
there holds $C = C^\alpha C^\beta$,
and $C^\alpha$ has,
since $y$ is relatively prime to $a$,
the order $a$,
and $C^\beta$ the order $b$.
(\textsc{Cauchy}, loc.~cit., \S.~V, p.~179.)
\origpagebreak
Let now also $C = A B$,
where $A$ and $B$ have the orders $a$ and $b$
and are permutable with each other.
Then
$C^\alpha = A^\alpha B^\alpha$,
$B^\alpha = B^{b y} = E$,
$A^\alpha = A^{1-\beta} = A$,
thus
$A = C^\alpha$ and $B = C^\beta$.
Being powers of $C$,
$A$ and $B$ belong to every group
to which $C$ belongs.

II.
\phantomsection\label{t:2-2}
\emph{ %
If the order of a group is divisible by $n$
then the number of those elements of the group
whose order divides $n$
is a multiple of $n$.
}

Let $\CH$ be a group of order $h$
and
$n$ a divisor of $h$.
For every group whose order is $h' < h$
and
for each divisor $n'$ of $h'$,
I assume the Theorem as proven.
The number of elements of $\CH$
whose order divides $n$
is,
if $n = h$ holds,
equal to $n$.
So if $n < h$,
I can assume
the theorem has been proven
for every divisor of $h$
which is $> n$.
Now if $p$ is a prime
dividing
$\frac{h}{n}$,
then
the number of elements of $h$
whose order divides $n p$
is divisible by $n p$,
hence also by $n$.
Let $n p = p^\lambda r$,
where $r$ is not divisible by $p$ and $\lambda \geq 1$.
Let $\CK$
be the complex
of those elements of $\CH$
whose order divides $n p$
but not $n$,
hence divisible by $p^\lambda$,
and let $k$ be the order of this complex.
Then it only remains to show
that the number $k$,
if it differs from zero,
is divisible by $n$.
For that purpose I prove
that $k$ is divisible
by $p^{\lambda-1}$ and $r$.

I partition the elements of $\CK$
into systems
by
assigning two elements
to the same system
if each is a power of the other.
All elements of a system
have the same order $m$.
Their number is $\phi(m)$.
A system is completely determined
by each of its elements $A$,
it is formed by the elements $A^\mu$
where $\mu$ runs through the $\phi(m)$ numbers
which are $< m$ and relatively prime to $m$.
If $A$ is an element of the complex $\CK$
then all the elements
of the system represented by $A$
belong to the complex $\CK$.
Then the order $m$ of $A$
is divisible by $p^\lambda$,
hence also $\phi(m)$ by $p^{\lambda-1}$.
Since the number of elements
of each system,
into which $\CK$ is decomposed,
is divisible by $p^{\lambda-1}$,
so
must $k$ be divisible by $p^{\lambda-1}$.

To show secondly
that
$k$ is also divisible by $r$,
I partition again
the elements of $\CK$ into systems,
but of a different kind,
yet still such that
the cardinality of elements
of each system is divisible by $r$.
Every element of $\CK$ can,
and in a unique way at that,
be represented
as a product
of an element $P$ of order $p^\lambda$
and
a with it permutable element $Q$
whose order divides $r$.
\origpagebreak
Conversely,
every product $P Q$
so obtained
belongs
to the complex $\CK$.

Let $P$ be some element of order $p^\lambda$.
All elements of $\CK$ that are permutable with $P$
form a group $\CQ$
whose order $q$ is divisible by $p^\lambda$.
The powers of $P$
form a group $\CP$ of order $p^\lambda$
which is an invariant subgroup of $\CQ$.
The elements $Q$ of $\CQ$
that satisfy the equation $Y^r = E$
are identical
to those
that satisfy the equation $Y^t = E$,
where $t$
is
the greatest common divisor of $q$ and $r$.
The first issue
is to determine
the number of those elements.

Every element of $\CQ$ can,
and in a unique way at that,
be represented
as a product
of an element $A$
whose order is a power of $p$
and
a with it permutable element $B$
whose order is not divisible by $p$.

If the $t$-th power of $A B$
belongs to the group $\CP$
then
\[
	(A B)^t = A^t B^t = P^s,
	\quad
	\text{hence}
	\quad
	A^t = P^s,
	\quad
	B^t = E,
\]
because also this element
can be decomposed in the given fashion
in a single way.
Thus
$A^t$ belongs to $\CP$,
hence also $A$ itself
because $t$ is not divisible by $p$.
The order of the group $\CQ / \CP$
is
$\frac{q}{p^\lambda} < h$.
The number of (complex) elements
of this group
that satisfy the equation $Y^t = R$
is therefore
a multiple of $t$,
say $t u$.
If $\CP A B$ is such an element
then,
as $A$ belongs to $\CP$,
$\CP A = \CP$,
hence $\CP A B = \CP B$.
Since $B$,
as an element of $\CQ$,
is permutable with $P$,
the complex $\CP B$
contains only one element
whose order divides $t$,
namely $B$ itself,
whilst
the order of every other element of $\CP B$
is divisible by $p$.
Let
\[
	\CP \CB + \CP B_1 + \CP B_2 + \cdots
\]
be the $t u$ distinct (complex) elements
of the group $\CQ / \CP$
whose
$t$-th power is contained in $\CP$,
then
this complex contains all those elements of $\CQ$
whose $t$-th power
(absolutely)
equals $E$.
However,
only the elements $B$, $B_1$, $B_2$, $\cdots$
have this property.
Thus $\CQ$ contains
exactly $t u$ elements
that satisfy the equation $Y^t = E$,
or
there are,
if $P$ is a certain element of order $p^\lambda$,
exactly $t u$ elements
that are permutable with $P$
and
whose order divides $r$.

The number of elements of $\CH$
permutable with $P$ is $q$.
The number of elements $P, P_1, P_2, \cdots$
of $\CH$
that are conjugate to $P$
with respect to $\CH$
is therefore $\frac{h}{q}$.
Then there are exactly $t u$
elements $Q_1$ in $\CH$
that are permutable with $P_1$
and
whose order divides $r$.
\origpagebreak
Taking
each of the $\frac{h}{q}$ elements
$P, P_1, P_2, \cdots$
successively as $X$
and
each time
as $Y$
the $t u$ elements permutable with $X$
and
satisfy the equation $Y^r = E$,
one obtains
the system
$\CK'$
of 
\[
	k' = \frac{h}{q} \, t u
\]
distinct elements $X Y$ of the complex $\CK$.
Now $h$ is divisible
by both $q$ and $r$
hence
also by their
least common multiple $\frac{q r}{t}$.
Thus $k'$ is divisible by $r$.
The system $\CK'$
is completely determined
by each of its elements.
Two distinct systems
among $\CK', \CK'', \cdots$
have no element in common.
Their orders $k', k'', \cdots$
are all divisible by $r$.
Thus also $k = k' + k'' + \cdots$ is
divisible by $r$.

The number of elements of a group
that satisfy the equation $X^n = E$
is $m n$,
the integer number $m$ is $> 0$
because
$X = E$ always satisfies that equation.

III.
\phantomsection\label{t:2-3}
\emph{ %
If the order of a group $\CH$ is divisible by $n$
then
the elements of $\CH$
whose order divides $n$
generate
a characteristic subgroup
of $\CH$
whose order is divisible by $n$.
}

Let $\CR$ be the complex of elements of $\CH$
that satisfy the equation $X^n = E$.
If $X$ is an element of $\CR$
and $R$ is any element%
\footnote{%
	\scriptsize
	tn: cf.~Frobenius,
	\emph{\"Uber endliche Gruppen}, \S.~5,
	SB.~Akad.~Berlin,
	1895 (I),
	\href{http://dx.doi.org/10.3931/e-rara-18846}{http://dx.doi.org/10.3931/e-rara-18846}
}
permutable
with $\CH$
then
$R^{-1} X R$
is also an element of $\CR$.
Thus $R^{-1} \CR R = \CR$.
Let the complex $\CR$ generate a group $\CG$
of order $g$.
Then also $R^{-1} \CG R = \CG$,
so that
$\CG$ is a characteristic subgroup of $\CH$.

If $q^\mu$ is the highest power
of a prime $q$
that divides $n$
then $q^\mu$ also divides $h$.
Thus $\CH$ contains a group $\CQ$ of order $q^\mu$.
Now $\CR$ is divisible by $\CQ$,
hence also $\CG$,
and consequently
$g$ is divisible by $q^\mu$.
Since this holds for every prime $q$ that divides $n$,
$g$ is divisible by $n$.

On the relation of the complex $\CR$ to the group $\CG$
I further note the following:
I considered
in \emph{\"Uber endliche Gruppen}, \S.~1
the powers $\CR, \CR^2, \CR^3, \cdots$ of a complex $\CR$.
If in that sequence
$\CR^{r + s}$
is the first one
that equals one of the foregoing ones $\CR^r$,
then
$\CR^\rho = \CR^\sigma$
if and only if
$\rho \equiv \sigma \pmod{s}$
and
$\rho$ and $\sigma$ are both $\geq r$.
Let $t$ be the number uniquely defined by
the conditions
$t \equiv 0 \pmod{s}$
and
$r \leq t < r + s$.
Then $\CR^t$ is the only group contained in that sequence
of powers.
If $\CR$ contains the principal element $E$
then $\CR^{\rho + 1}$ is divisible by $\CR^\rho$.
Hence $\CG = \CR^t$ is divisible by $\CR$.
If $N$ is an element of the group $\CG$
then $\CG N = \CG$.
More generally then,
if $\CR$
is a complex of elements
contained in $\CG$
then $\CG \CR = \CG$.
Therefore $\CR^{t+1} = \CR^t$,
hence $s = 1$ and $t = r$.
\origpagebreak
Consequently,
$\CR^r = \CR^{r+1}$
is the first one
in the sequence of powers of $\CR$
that equals the subsequent one,
and this is the group generated by the complex $\CR$.

IV.
\phantomsection\label{t:2-4}
\emph{ %
If the order of a group $\CH$
is
divisible by the two relatively prime numbers $r$ and $s$,
if
there exists in $\CH$
exactly $r$ elements $A$
whose order divides $r$
and
exactly $s$ elements $B$
whose order divides $s$,
then
each of the $r$ elements $A$
is permutable
with each of the $s$ elements $B$
and
there exist in $\CH$ exactly $r s$ elements
whose order divides $r s$,
namely the $r s$ distinct elements $A B = B A$.
}

Indeed,
every element $C$ of $\CH$
whose order divides $r s$
can be written as a product of
two with each other permutable elements $A$ and $B$
whose orders divide $r$ and $s$.
Now
$\CH$ contains
no more than $r$ elements $A$
and
no more than $s$ elements $B$.
Were it not the case
that
each of the $r$ elements $A$
is permutable with
each of the $s$ elements $B$
and
furthermore
that
the $r s$ elements $A B$ are all distinct,
then
$\CH$ would contain less than $r s$ elements $C$.
But this contradicts
\hyperref[t:2-2]{Theorem II}.

\subsubsection*{\secformat{\S.~3.}}

If the order $h$ of a group $\CH$ divisible by the prime $p$
then
$\CH$ contains elements of order $p$,
namely $m p - 1$ many,
because there exist $m p$ elements in $\CH$
whose order divides $p$.
From this theorem of \textsc{Cauchy},
\textsc{Sylow}
derived the more general one,
that
any group whose order is divisible by $p^\kappa$
possesses
a subgroup of order $p^\kappa$.
In his proof he draws on
the language of the theory of substitutions.
If one wants to avoid this,
one should apply the procedure
that I used
in my work
\emph{\"Uber endliche Gruppen}
in the proof of
Theorems V and VII, \S.~2.

Another proof is obtained
by partitioning
the $m p - 1$ elements $P$ of order $p$ contained in $\CH$
into classes of conjugate elements.
If the elements of $\CH$ permutable with $P$
form a group $\CG$ of order $g$,
then
the number of elements conjugate to $P$ is~$\frac{h}{g}$.
Thus
\[
	m p - 1
	=
	\sum \frac{h}{g}
\]
where the sum is to be extended
over the different classes
into which the elements $P$ are segregated.
From this equation it follows
that
not all the summands $\frac{h}{g}$
are divisible by $p$.
Let $p^\lambda$ be the highest power of $p$ contained in $h$,
and
let $\kappa \leq \lambda$.
If $\frac{h}{g}$ is not divisible by $p$
then $g$ is divisible by $p^\lambda$.
\origpagebreak
The powers of $P$ form a group $\CP$ of order $p$,
which is an invariant subgroup of $\CG$.
The order of the group $\CG / \CP$ is $\frac{g}{p} < h$.
For this group we may therefore
assume the theorems
which we wish to prove for $\CH$ as known.
Thus it contains
a group $\CP_\kappa / \CP$
of order $p^{\kappa - 1}$,
and in the case that $\kappa < \lambda$,
a group $\CP_{\kappa+1} / \CP$ of order $p^\kappa$
that is divisible by $\CP_\kappa / \CP$.
Consequently,
$\CH$ contains the group $\CP_\kappa$ of order $p^\kappa$
and the group $\CP_{\kappa+1}$ of order $p^{\kappa+1}$
that is divisible by $\CP_\kappa$.

\subsubsection*{\secformat{\S.~4.}}

I.
\phantomsection\label{t:4-1}
\emph{ %
If the order of a group is divisible by
the $\kappa$-th power of the prime $p$
then
the number of groups of order $p^\kappa$
contained therein
is a number of the form $n p + 1$.
}

Let $r_\kappa$ denote the number of groups of order $p^\kappa$
contained in $\CH$.
Then the number of elements of $\CH$
whose order is $p$
equals $r_1 (p - 1)$.
As shown above,
this number has the form $m p - 1$.
Thus
\begin{align} \label{e:4-1} \tag{1.}
	r_1 \equiv 1 \pmod{p}
	.
\end{align}

Let $r_{\kappa - 1} = r$, $r_\kappa = s$,
and let
\begin{align} \label{e:4-2} \tag{2.}
	\CA_1, \CA_2, \cdots, \CA_r
\end{align}
be the $r$ groups of order $p^{\kappa-1}$ contained in $\CH$
and
\begin{align} \label{e:4-3} \tag{3.}
	\CB_1, \CB_2, \cdots, \CB_s
\end{align}
the $s$ groups of order $p^\kappa$.
Suppose the group $\CA_\rho$
is contained in
$a_\rho$ of the groups \eqref{e:4-3}.
Suppose the group $\CB_\sigma$
is divisible by
$b_\sigma$ of the groups \eqref{e:4-2}.
Then
\begin{align} \label{e:4-4} \tag{4.}
	a_1 + a_2 + \cdots + a_r
	=
	b_1 + b_2 + \cdots + b_s
\end{align}
is the number of distinct pairs of groups
$\CA_\rho, \CB_\sigma$
for which $\CA_\rho$ is contained in $\CB_\sigma$.

Let $\CA$ be one of the groups \eqref{e:4-2}.
Of the groups \eqref{e:4-3}
let
$\CB_1, \CB_2, \cdots, \CB_a$
be those
that are divisible by $\CA$.
By \S.~3,
$a > 0$,
and by
\hyperref[t:1-2]{Theorem II, \S.~1},
$\CA$ is an invariant subgroup of each of these $a$ groups,
hence also of their
least common multiple $\CG$.
Therefore
the group $\CG / \CA$
contains the $a$ groups
$\CB_1 / \CA, \CB_2 / \CA, \cdots, \CB_a / \CA$
of order $p$
and none further.
Indeed,
if $\CB / \CA$ is a group of order $p$
contained in $\CG / \CA$
then
$\CB$ is a group of order $p^\kappa$
divisible by $\CA$.
\origpagebreak
By formula \eqref{e:4-1}
there holds
$a \equiv 1 \pmod{p}$.
Thus
\begin{align} \label{e:4-5} \tag{5.}
	a_\rho \equiv 1,
	\quad
	a_1 + a_2 + \cdots + a_r \equiv r
	\pmod{p}
	.
\end{align}

Now I need the Lemma:

\emph{The number of groups of order $p^{\lambda-1}$
which are contained in a group of order $p^\lambda$
is $\equiv 1 \pmod{p}$.}

I suppose this Lemma is already proven
for groups of order $p^\kappa$ if $\kappa < \lambda$.
Then,
if in the above expansion $\kappa < \lambda$
then
\begin{align} \label{e:4-6} \tag{6.}
	b_\sigma \equiv 1,
	\quad
	b_1 + b_2 + \cdots + b_s \equiv s
	\pmod{p}.
\end{align}

Therefore $r \equiv s$ or $r_{\kappa-1} \equiv r_\kappa \pmod{p}$,
and since this congruence
holds for each value $\kappa < \lambda$,
it is
\[
	1 \equiv r_1 \equiv r_2 \equiv \cdots \equiv r_{\lambda-1}
	\pmod{p}
	.
\]

Applying this result
to a group $\CH$
whose order is $p^\lambda$,
it is therefore
$r_{\lambda-1} \equiv 1 \pmod{p}$
for such a group,
and
with this,
the above Lemma is proven
also for groups of order $p^\lambda$,
if it holds
for groups of order $p^\kappa < p^\lambda$,
it is therefore generally valid.
For each value $\kappa$
consequently,
$r_\kappa \equiv r_{\kappa-1}$
and
therefore
$r_\kappa \equiv 1 \pmod{p}$.

In exactly the same way
one proves the more general Theorem:

II.
\phantomsection\label{t:4-2}
\emph{ %
If the order of a group $\CH$
is divisible by the $\kappa$-th power
of the prime $p$,
if $\vartheta \leq \kappa$
and
$\CP$ is a group of order $p^\vartheta$
contained in $\CH$,
then
the number of groups of order $p^\kappa$
contained in $\CH$
that are divisible by $\CP$
is a number of the form $n p + 1$.
}

\subsubsection*{\secformat{\S.~5.}}

The Lemma used in \S.~4
can also be proven
in the following way
by relying on the Theorem:
Every group $\CH$ of order $p^\lambda$
has a subgroup $\CA$ of order $p^{\lambda-1}$
and such a subgroup
is always an invariant one.
Let $\CA$ and $\CB$ be
two distinct subgroups of order $p^{\lambda-1}$
contained in $\CH$
and
let $\CD$ be their greatest common divisor.
Since $\CA$ and $\CB$ are invariant subgroups of $\CH$,
so is $\CD$,
and
since $\CH$ is the least common multiple of $\CA$ and $\CB$,
$\CD$ has order $p^{\lambda-2}$.
Thus
$\CH / \CD$ is a group of order $p^2$.
Any such group has,
depending on whether
it is a cyclic group or not,
$1$ or $p+1$
subgroups of order $p$,
thus in our case
$p + 1$,
since $\CA / \CD$ and $\CB / \CD$
are two distinct groups of this type.
Therefore
$\CH$ contains exactly $p + 1$
distinct groups of order $p^{\lambda-1}$
that are divisible by $\CD$.

\origpagebreak

The group $\CH$ always contains
a group $\CA$ of order $p^{\lambda-1}$.
If it contains yet another one,
then
$\CH$ has an invariant subgroup $\CD$ of order $p^{\lambda-2}$
which is contained in $\CA$
and
for which the group $\CH / \CD$ is not a cyclic one.
Let $\CD_1, \CD_2, \cdots, \CD_n$
be
all the groups of this kind.
Then there exist in $\CH$ besides $\CA$
other $p$ groups of order $p^{\lambda-1}$
divisible by $\CD_1$
\begin{align} \label{e:5-1} \tag{1.}
	\CA_1, \CA_2, \cdots, \CA_p
	,
\end{align}
and likewise
$p$ groups that are
divisible by $\CD_2$
\begin{align} \label{e:5-2} \tag{2.}
	\CA_{p+1}, \CA_{p+2}, \cdots, \CA_{2 p}
	,
\end{align}
etc.,
and finally
$p$ groups
divisible by $\CD_n$
\begin{align} \label{e:5-3} \tag{3.}
	\CA_{(n-1) p + 1},
	\CA_{(n-1) p + 2},
	\cdots,
	\CA_{n p + 1}
	.
\end{align}

The $n p + 1$ groups
$\CA, \CA_1, \cdots, \CA_{n p}$
are all the groups of order $p^{\lambda-1}$ contained in $\CH$
since
each such group $\CB$
has to have
in common with $\CA$
a certain divisor $\CD$
which is one of the $n$ groups $\CD_1, \CD_2, \cdots, \CD_n$.
They are, furthermore, all distinct.
Indeed,
if $\CA_1 = \CA_{p+1}$ was true
then
$\CA_1$ would be divisible by both groups $\CD_1$ and $\CD_2$,
hence also by their least common multiple $\CA$.
If $\CP$ is a group of order $p^\vartheta$
contained in $\CH$
then
one can subject all the groups considered above
to the condition
of being divisible by $\CP$.
If conversely
$\CH$ is an invariant subgroup of a group $\CP$ of order $p^\vartheta$
then
one can require
that
they all be invariant subgroups of $\CP$.

With the help of
\hyperref[t:1-5]{Theorem V, \S.1}
it is easy to prove
that the number of groups of order $p^{\lambda-1}$
that are contained in a group of order $p^\lambda$
equals $1$
only if
$\CH$ is a cyclic group.

I.
\phantomsection\label{t:5-1}
\emph{ %
The number of invariant subgroups of order $p^\kappa$
contained in a group of order $p^\lambda$
is a number of the form $n p + 1$.
}

Let $\CH$ be a group of order $h$,
let $p^\lambda$ be the highest power of $p$ contained in $h$,
let $\kappa \leq \lambda$
and
$\CP_\kappa$ any group of order $p^\kappa$
contained in $\CH$.
Each group $\CP_\kappa$ is contained in $n p + 1$ groups,
hence at least in one.
I divide the groups $\CP_\kappa$
into two kinds.
For a group of the first kind
there exists
a group $\CP_\lambda$
of which $\CP_\kappa$ is an invariant subgroup,
for a group of the second kind
no such group exists.
The number of elements of $\CH$
permutable with $\CP_\kappa$
is divisible by $p^\lambda$ in the first case,
and
in the second case it is not.
The number of groups
conjugate to $\CP_\kappa$ is therefore
divisible by $p$ in the second case,
in the first case it is not.
Hence diving the groups $\CP_\kappa$
into classes of conjugate groups
one recognizes
that the number of groups $\CP_\kappa$
of the second kind
is divisible by $p$.
\origpagebreak
Consequently,
the number of groups $\CP_\kappa$ of the first kind
is $\equiv 1 \pmod{p}$.

II.
\phantomsection\label{t:5-2}
\emph{ %
If $\CH$ is a group of order $p^\lambda$
and $\CG$ is an invariant subgroup of $\CH$
whose order is divisible by $p^\kappa$
then the number of groups of order $p^\kappa$
contained in $\CG$
that are invariant subgroups of $\CH$
is
a number of the form $n p + 1$.
}

Also here let more generally
$p^\lambda$ be the highest power of the prime $p$
that divides the order $h$ of $\CH$.
Let $\CG$ be an invariant subgroup of $\CH$
whose order $g$ is divisible by $p^\kappa$.
The number of all groups $\CP_\kappa$ of order $p^\kappa$
contained in $\CG$
is $\equiv 1 \pmod{p}$.
I divide them into groups of the first and the second kind
(with respect to $\CH$)
and
further
into classes of conjugate groups.
If $\CG$ is divisible by $\CP_\kappa$
then
$\CG$ is also divisible by every group
conjugate to $\CP_\kappa$.
Therefrom
the claim follows in the same way as above.
One can also easily prove it directly
by means of the method
used in \S.~4:

Let the order of $\CH$ be $h = p^\lambda$.
By
\hyperref[t:1-5]{Theorem V, \S.1}
the group $\CG$ contains elements of order $p$
that are invariant elements of $\CH$.
They form, together with the principal element,
a group.
If $p^\alpha$ is its order
then $p^\alpha - 1$ is the number of those elements.
By
\hyperref[t:1-3]{Theorem III, \S.1},
every invariant subgroup of $\CH$
whose order is $p$
consists of
the powers of such an element.
Therefore
there exist in $\CG$
$r = \frac{p^\alpha - 1}{p - 1}$
groups of order $p$
that are invariant subgroups of $\CH$.
This number is
\begin{align} \label{e:5-4} \tag{4.}
	r \equiv 1 \pmod{p}
	.
\end{align}
Let
\begin{align} \label{e:5-5} \tag{5.}
	\CA_1, \CA_2, \cdots, \CA_r
\end{align}
be those $r$ groups
and let
\begin{align} \label{e:5-6} \tag{6.}
	\CB_1, \CB_2, \cdots, \CB_s
\end{align}
be the $s$ groups of order $p^\kappa$
contained in $\CG$
that are invariant subgroups of $\CH$.
Let $\CB$ be one of the groups \eqref{e:5-6}.
Among the groups \eqref{e:5-5}
let
$\CA_1, \CA_2, \cdots, \CA_b$
be those contained in $\CB$.
By \eqref{e:5-4}
is then
$b \equiv 1 \pmod{p}$.
Let $\CA$ be one of the groups \eqref{e:5-5}.
Among the groups \eqref{e:5-6}
let $\CB_1, \CB_2, \cdots, \CB_a$
be those divisible by $\CA$.
Then
$\CB_1 / \CA, \CB_2 / \CA, \cdots, \CB_a / \CA$
are
the groups of order $p^{\kappa-1}$
contained in $\CG / \CA$
that are invariant subgroups of $\CH / \CA$.
By the method of induction
is therefore
$a \equiv 1 \pmod{p}$.
Resorting to the same notation as in \S.~4
there holds
\[
	1 \equiv r
	\equiv
	a_1 + a_2 + \cdots + a_r
	\equiv
	b_1 + b_2 + \cdots + b_s
	\equiv
	s
	\pmod{p}
	.
\]

\origpagebreak

I add a few remarks
on the number of groups $\CP_\kappa$
of the first kind
that
are conjugate to a particular one,
and
on the number of classes of conjugate groups
into which
the groups $\CP_\kappa$
are partitioned.

Let $\CP$ be a group of order $p^\lambda$
contained in $\CP$
and
$\CQ$ an invariant subgroup of $\CP$ of order $p^\kappa$.
The elements of $\CH$ permutable with $\CP$ $(\CQ)$
form a group of $\CP'$ $(\CQ')$
of order $p'$ $(q')$.
Let
the greatest common divisor
of $\CP'$ and $\CQ'$
be the group $\CR$ of order $r$.
The groups
$\CP'$, $\CQ'$ and $\CR$
are divisible by $\CP$.
Let $p^\delta$
be the order of
the largest common divisor
of $\CP$
and
a group conjugate
with respect to $\CH$
that is selected
in such a way
that $\delta$ is a maximum.
Then
(\emph{\"Uber endliche Gruppen}, \S.~2, VIII)
\[
	\frac{h}{p'}
	\equiv
	1
	\pmod{p^{\lambda-\delta}}
	.
\]
The group $\CR$
consists
of all the elements of $\CQ'$
that are permutable with $\CP$.
With this,
\[
	\frac{q'}{r} \equiv 1
	\pmod{p^{\lambda-\delta}}
	.
\]
Consequently,
\begin{align} \label{e:5-7} \tag{7.}
	\frac{h}{q'} \equiv \frac{p'}{r}
	\pmod{p^{\lambda-\delta}}
	.
\end{align}

Herein,
$\frac{h}{q'}$
is the number of groups
that are conjugate to $\CQ$
with respect to $\CH$
and
$\frac{p'}{r}$
is the number of groups
that are conjugate to $\CQ$
with respect to $\CP'$.
Indeed,
the group $\CR$
consists
of all the elements of $\CP'$
that are permutable with $\CQ$.
The number of groups
in a certain class in $\CH$
is therefore
congruent $\!\!\pmod{p^{\lambda-\delta}}$
to
the number of groups
in the corresponding class in $\CP'$.

Furthermore,
the number of distinct classes in $\CH$
(into which the groups $\CP_\kappa$ of the first kind
are partitioned)
equals
the number of those classes in $\CP'$.
This follows from the Theorem:

III.
\phantomsection\label{t:5-3}
\emph{ %
If two invariant subgroups of $\CP$
are conjugate with respect to $\CH$
then
so they are with respect to $\CP'$.
}

Let $\CQ$ and $\CQ_0$
be two invariant subgroups of $\CP$.
If they are conjugate with respect to $\CH$
then
there exists in $\CH$
such an element $H$
that
\begin{align} \label{e:5-4'} \tag{4.}
	H^{-1} \CQ_0 H = \CQ
\end{align}
holds.
Since $\CQ_0$ is an invariant subgroup of $\CP$,
$H^{-1} \CQ_0 H = \CQ$
is an invariant subgroup of
\[
	H^{-1} \CP H = \CP_0
	.
\]
\origpagebreak
Hence
$\CQ'$ is divisible by $\CP$ and $\CP_0$.
Consequently
(\emph{\"Uber endliche Gruppen}, \S.~2, VII)
there exists in $\CQ'$
such an element $Q$
that
\[
	Q^{-1} \CP_0 Q = \CP
	,
\]
hence
\[
	\CP H Q = H Q \CP
\]
holds.
Thus $H Q = P$
is an element of $\CP'$.
Inserting the expression $H = P Q^{-1}$
into the equation \eqref{e:5-4'}
one obtains,
since $Q$ is permutable with $\CQ$,
\[
	P^{-1} \CQ_0 P
	=
	Q^{-1} \CQ Q
	=
	\CQ
	.
\]
There exists therefore
in $\CP'$
an element $P$
that transforms $\CQ_0$ into $\CQ$.

Partition now
the groups $\CP_\kappa$
contained in $\CH$
(of the first kind)
into classes
of conjugate groups
(with respect to $\CH$)
and
choose
from each class a representative.
If $\CQ_0$ is one,
then $\CQ_0$ is a group of order $p^\kappa$
which
is contained in a certain group $\CP_0$
as an invariant subgroup.
If $H^{-1} \CP_0 H = \CP$
then $H^{-1} \CQ_0 H = \CQ$
is an invariant subgroup of $\CP$.
One can therefore
choose the representatives
of different classes
in such a way
that
they are all invariant subgroups
of a certain group $\CP$ of order $p^\lambda$.
Each invariant subgroup of $\CP$
of order $p^\kappa$
is then conjugate
to one of these groups
with respect to $\CH$,
hence also
with respect to $\CP'$.
Let
the invariant subgroups $\CP_\kappa$ of $\CP$
aggregate
into $s$ classes of groups
that are conjugate
with respect to $\CP'$.
Then the groups $\CP_\kappa$
of the first kind
of $\CH$
also
aggregate
into $s$ classes of groups
that are conjugate
with respect to $\CH$.

\vfill

\noindent
\begin{minipage}{\textwidth}
{\footnotesize

\noindent
Translated with minor typographical corrections from
\begin{quote}
	F.~G.~Frobenius,
	``Verallgemeinerung des Sylow'schen Satzes'',
	Sitzungsberichte der
	K\"onigl.~Preu\ss.~Akad.~der Wissenschaften
	zu Berlin,
	1895 (II),
	981--993.
	\\
	\href{http://dx.doi.org/10.3931/e-rara-18880}{http://dx.doi.org/10.3931/e-rara-18880}
\end{quote}

\noindent
For the terminology see
\begin{quote}Frobenius,
	\emph{\"Uber endliche Gruppen},
	SB.~Akad.~Berlin,
	1895 (I),
	163--194.
	\\
	\href{http://dx.doi.org/10.3931/e-rara-18846}{http://dx.doi.org/10.3931/e-rara-18846}
\end{quote}

R.~Andreev,
Paris (France),
\today.
}
\end{minipage}


\end{document}